\documentclass[12pt]{article}
\RequirePackage{algorithm,algorithmic,amssymb,epsf,epic}

\def\denseformat{
\setlength{\textheight}{9in}
\setlength{\textwidth}{6.9in}
\setlength{\evensidemargin}{-0.2in}
\setlength{\oddsidemargin}{-0.2in}
\setlength{\headsep}{10pt}
\setlength{\topmargin}{-0.3in}
\setlength{\columnsep}{0.375in}
\setlength{\itemsep}{0pt}
}

\newtheorem{theorem}{Theorem}[section]

\newtheorem{lemma}[theorem]{Lemma}

\newtheorem{comment}[theorem]{Comment}

\def\boldhead#1:{\par\vskip 7pt\noindent{\bf #1:}\hskip 10pt}
\def\ithead#1:{\par\vskip 7pt\noindent{\it #1:}\hskip 10pt}

\def\inline#1:{\par\vskip 7pt\noindent{\bf #1:}\hskip 10pt}
\def\midinline#1:{\par\noindent{\bf #1:}\hskip 10pt}
\def\dnsinline#1:{\par\vskip -7pt\noindent{\bf #1:}\hskip 10pt}
\def\ddnsinline#1:{\newline{\bf #1:}\hskip 10pt}
\def\largeinline#1:{\par\vskip 7pt\noindent{\large\bf #1:}\hskip 10pt}
\long\def\comment #1\commentend{}
\long\def\commhide #1\commhideend{}
\long\def\commfull #1\commend{#1}
\long\def\commabs #1\commenda{}
\long\def\commtim #1\commendt{#1}
\long\def\commb #1\commbend{}
\long\def\commedit #1\commeditend{} 

\long\def\commB #1\commBend{}       

\long\def\commex #1\commexend{}     

\long\def\commsiena #1\commsienaend{}  
                                         
\long\def\commBI #1\commBIend{}  
                                         
\long\def\CProof #1\CQED{}

\def\blackslug{\hbox{\hskip 1pt \vrule width 4pt height 8pt
    depth 1.5pt \hskip 1pt}}
\def\QED{\quad\blackslug\lower 8.5pt\null\par}
\def\inQED{\quad\quad\blackslug}

\def\Proof{\par\noindent{\bf Proof:~}}

\def\proof{\Proof}

\long\def\PPP#1{\noindent{\bf Proof:}{ #1}{\quad\blackslug\lower 8.5pt\null}}

\long\def\denspar #1\densend
{#1}

\newif\ifnotesw\noteswtrue
   {\ifnotesw\marginpar[\hfill\(\top\)]{\(\top\)}\fi}%
   {\ifnotesw\marginpar[\hfill\(\bot\)]{\(\bot\)}\fi}

\newcommand{\mnote}[1]%
    {\ifnotesw\marginpar%
        [{\scriptsize\it\begin{minipage}[t]{\marginparwidth}
        \raggedleft#1%
                        \end{minipage}}]%
        {\scriptsize\it\begin{minipage}[t]{\marginparwidth}
        \raggedright#1%
                        \end{minipage}}%
    \fi}

\def\cB{{\cal B}}
\def\cC{{\cal C}}

\def\cH{{\cal H}}

\def\cS{{\cal S}}

\def\hD{{\hat D}}

\def\hN{{\hat N}}

\def\hW{{\hat W}}

\def\hZ{{\hat Z}}

\def\hgamma{{\hat \gamma}}

\def\bA{{\bar A}}

\def\bV{{\bar V}}

\def\tO{{\tilde O}}

\def\tQ{{\tilde Q}}
\def\tR{{\tilde R}}

\def\tW{{\tilde W}}

\def\tZ{{\tilde Z}}
\def\dist{{\mathit{dist}}}
\def\MathF{\hbox{\rm I\kern-2pt F}}
\def\MathP{\hbox{\rm I\kern-2pt P}}
\def\MathR{\hbox{\rm I\kern-2pt R}}
\def\MathZ{\hbox{\sf Z\kern-4pt Z}}
\def\MathN{\hbox{\rm I\kern-2pt I\kern-3.1pt N}}
\def\MathC{\hbox{\rm \kern0.7pt\raise0.8pt\hbox{\footnotesize I}
\kern-4.2pt C}}
\def\MathQ{\hbox{\rm I\kern-6pt Q}}

\def\MathE{\hbox{{\rm I}\hskip -2pt {\rm E}}} 

\newsavebox{\ttop}\newsavebox{\bbot}

\def\mod{\pmod}

\def\eps{\epsilon}

\def\setmns{\setminus}

\newcommand{\half}{{1\over2}}

\def\nin{{~\not \in~}}
\def\emset{\emptyset}

\newcommand{\Prob}{\MathP}
\newcommand{\Expect}{\MathE}

\def\etal{\emph{et~al.}}
\def\Erdos{Erd\H{o}s}  
\def\Szemeredi{Szemer\'edi}
\def\Turan{Tur\'an}
\denseformat

\begin{document}

\def\tO{\tilde{O}}
\def\nin{\not \in}
\def\emset{\emptyset}
\def\setmns{\setminus}
\def\etal{{et al.~}}
\def\Sp{\mathit{Sp}}
\def\ReadDg{\mathit{Read\_Edge}}
\def\wmax{{{\hat \omega}}}
\def\ttl{\mathit{ttl}}
\def\Rnd{\mathit{Round}}
\def\deg{\mathit{deg}}
\def\ctr{\mathit{ctr}}
\def\ScDgs{\mathit{Scanned\_Edges}}
\def\SCD{\mathit{SCANNED}}
\def\NOTSCD{\mathit{NOTSCANNED}}
\def\CRASH{\mathit{CRASH}}
\def\CR{\mathit{CRASH}}
\def\CRASHED{\mathit{CRASHED}}
\def\SyncIncr{\mathit{Sync\_Incr}}
\def\actdeg{\mathit{actdeg}}
\def\OLD{\mathit{OLD}}
\def\NEW{\mathit{NEW}}
\def\mark{\mathit{mark}}
\def\AsyncRnd{\mathit{Async\_Rnd}}
\def\status{\mathit{status}}
\def\scstat{\mathit{scan\_status}}
\def\lab{\mathit{label}}
\def\seclabel{\mathit{sec\_label}}
\def\own{\mathit{own}}
\def\SELF{\mathit{SELF}}
\def\PEER{\mathit{PEER}}
\def\DynRnd{\mathit{Dyn\_Rnd}}
\def\Delete{\mathit{Delete}}
\def\XReplace{\mathit{XReplace}}
\def\UpdLab{\mathit{Update\_Label}}
\def\Crash{\mathit{Crash}}
\def\CrashLoop{\mathit{CrashLoop}}
\def\CrashItern{\mathit{CrashItern}}
\def\fetched{\mathit{fetched}}
\def\TRUE{\mathit{TRUE}}
\def\FALSE{\mathit{FALSE}}
\def\crashstat{\mathit{crash\_status}}
\def\wmax{\hat{\omega}}
\def\END{\mathit{END}}
\def\te{\tilde{e}}
\def\tu{\tilde{u}}
\def\che{\check{e}}
\def\chu{\check{u}}
\def\Old{\mathit{Old}}
\def\tR{\tilde{R}}
\def\eps{{\epsilon}}
\def\UpLab{\mathit{Update\_Label}}
\def\chP{\check{P}}
\def\VOID{\mathit{VOID}}
\def\dist{\mathit{dist}}
\def\mod{\pmod}

\def\Var{\mathit{Var}}
\def\Ext{\mathit{Ext}}
\def\Vl{\mathit{Vol}}

\def\hcS{\hat{\cal S}}
\def\tcS{\tilde{\cal S}}
\def\ccS{\check{\cal S}}
\def\chS{\check{S}}

\title{An Improved Construction of Progression-Free Sets
}
\author{
Michael Elkin
\thanks{Department of Computer Science, Ben-Gurion University of the
  Negev, Beer-Sheva, Israel, {\tt elkinm@cs.bgu.ac.il}
\newline
This research has been supported by the Israeli Academy of Science,
grant 483/06.}
}


\maketitle

\begin{abstract}
The problem of constructing dense subsets $S$ of $\{1,2,\ldots,n\}$
that contain no arithmetic triple was introduced by \Erdos 
$~$
 and \Turan $~$ in 1936. They have presented a construction with $|S| =
 \Omega(n^{\log_3 2})$ elements. Their construction was improved by
 Salem and Spencer, and further improved by Behrend in 1946.
The lower bound of Behrend is 
$$|S| ~=~ \Omega\left({ n \over {{2^{2 \sqrt{2} \sqrt{\log_2 n}}} \cdot
  \log^{1/4} n}} \right)~.$$
Since then the problem became one of the most central, most
  fundamental, and most intensively studied problems in additive
  number theory. Nevertheless, no improvement of the lower bound of
  Behrend was reported since 1946.

In this paper we present a construction that improves the result of
Behrend by a factor of $\Theta(\sqrt{\log n})$, and shows that
$$|S| ~=~ \Omega\left({ n \over {{2^{2 \sqrt{2} \sqrt{\log_2 n}}} }}
  \cdot \log^{1/4} n \right)~.$$
In particular, our result implies that the construction of Behrend is
  not optimal. 

Our construction is elementary and self-contained.

\end{abstract}



\section{Introduction}

A subset $S \subseteq \{1,2,\ldots,n\}$ is called {\em
  progression-free} if it contains no three distinct elements
  $i,j,\ell \in S$ such that $i$ is the arithmetic average of $j$ and
  $\ell$, i.e., $ i = {{j + \ell} \over 2}$. For a positive integer
  $n$, let $\nu(n)$ denote the largest size of a progression-free
  subset $S$ of $\{1,2,\ldots,n\}$.

Providing asymptotic estimates on
  $\nu(n)$ is a central and fundamental problem in additive number
  theory.
This problem was introduced by \Erdos $~$ and Turan \cite{ET36} in
  1936, and they have shown that $\nu(n) = \Omega(n^{\log_3 2})$.
This estimate was improved by Salem and Spencer \cite{SS42}, and
  further improved by Behrend \cite{B46} in 1946.
Behrend have shown that
$$\nu(n) ~=~ \Omega\left({n \over {{2^{2 \sqrt{2} \sqrt{\log n}}}
    \cdot \log^{1/4} n}}\right)~,$$
and this bound remains state-of-the-art for more than sixty years. A
    slightly weaker lower bound that does not rely on the Pigeonhole
    Principle was shown by Moser \cite{M53}.

The first non-trivial upper bound $\nu(n) = O({n \over {\log\log n}})$
was proved in a seminal paper by Roth \cite{R53}.
This bound was improved by Bourgain \cite{B99,B07}, and the current
state-of-the-art upper bound is 
$\nu(n) = O(n \cdot {{(\log \log n)^2} \over {\log^{2/3} n}})$
\cite{B07}.
The problem is also closely related to \Szemeredi $~$ theorem \cite{Sz75},
and to the problem of finding arbitrarily long arithmetic progressions of
prime numbers (see, e.g., Green and Tao \cite{GT07}), and to other
central problems in the additive number theory.
(See, e.g., the enlightening survey  on  \Szemeredi $~$ theorem in 
Scholarpedia.)

In this paper we improve the lower bound of Behrend by a factor of
$\Theta(\sqrt{\log n})$, and show that 
$$|S| ~=~ \Omega\left({ n \over {{2^{2 \sqrt{2} \sqrt{\log_2 n}}} }}
  \cdot \log^{1/4} n \right)~.$$
Though the improvement is not large, our result demonstrates that the
  construction of Behrend is not optimal. Also, despite very intensive
  research in this area, no improvement of Behrend lower bound was
  achieved for more than sixty years.

Our proof is elementary, and self-contained. The proof of Behrend is
based on the observation that a sphere in any dimension is convex, and
thus cannot contain an arithmetic progression. We replace the sphere
by a thin annulus, and demonstrate that this annulus contains a large
convexly independent subset $U$ of integer points.
There is an inherent tradeoff between the width of the annulus and the
size of $U$. In our construction we choose the largest width for which
we are able to show that $U$ contains at least a constant fraction of
all integer points of the annulus.

The construction of Behrend was generalized by Rankin \cite{Ran60} to
provide large subsets of $\{1,2,\ldots,n\}$ that contain no arithmetic
progression of length $k$, for any fixed $k$.
We believe that our technique will be useful for improving the lower
bound of Rankin as well.
Finally, like the construction of Behrend, our construction relies on
the Pigeonhole Principle. Consequently, the result of Moser \cite{M53}
remains the best known lower bound achieved without relying on the
Pigeonhole Principle. However, we hope that our argument can be made
independent of the Pigeonhole Principle. (See also Section
\ref{sec:concl}.) 

\section{Preliminaries}
\label{sec:prel}

For a pair $a,b$ of real numbers, $a \le b$, we denote by $[a,b]$
(respectively,  $(a,b)$)
the closed (resp.,  open) segment containing all numbers $x$, $a \le x
\le b$ (resp., $a < x < b$). We also use the notation $(a,b]$
  (respectively, 
  $[a,b)$) for denoting the segment containing all numbers $x$, $a < x
    \le b$ (resp., $a \le x < b$).
For  integer numbers $n$ and $m$, $n \le m$, we denote by $[\{n,m\}]$
the set of integer numbers $\{n,n+1,\ldots,m\}$. If $n = 1$ then we
use the notation $[\{m\}]$ as a shortcut for $[\{1,m\}]$.
For a real number $x$, we denote by $\lfloor x \rfloor$ (respectively,
$\lceil x \rceil$) the largest (resp., smallest)
integer number that is no greater (resp., no smaller) than $x$.

A triple $i,j,\ell$ of distinct integer numbers is called an {\em
  arithmetic triple}  if one of these numbers is the average of two
  other numbers, i.e., $i = {{j + \ell} \over 2}$.
A set $S$ of integer numbers is called {\em progression-free} if it
  contains no arithmetic triple. For a positive integer number $n$,
let $\nu(n)$ denote the largest size of a progression-free subset $S$
  of
$[\{n\}]$.

For a pair of integer functions $f(\cdot)$, $g(\cdot)$,  we say
that $f(n) = O(g(n))$ if there exists  a positive (universal)  constant
$c$ and a positive integer $N$   such that 
for every $n \ge N$, $|f(n)| \le
  c \cdot |g(n)|$. In this case we also say that $g(n) = \Omega(f(n))$.
If both $f(n) = O(g(n))$ and $g(n) = O(f(n))$ hold, we say that $f(n) =
\Theta(g(n))$.
If  $\lim_{n \to \infty} {{f(n)} \over
  {g(n)}} = 0$ we say that $f(n) = o(g(n))$.
These definitions extend to positive real functions as well.

Unless specified explicitly, $\log$ (respectively, $\ln$) stands for
the logarithm on {\em base 2} (resp., $e$).

For a positive integer $k$ and a  vector $v = (v_1,v_2,\ldots,v_k)$,
let $||v|| = \sqrt{\sum_{i=1}^k v_i^2}$ denote the {\em  norm} of
the vector $v$. The expression $||v||^2 = \sum_{i=1}^k v_i^2$ will be
referred to as the {\em squared norm} of the vector $v$.

For three vectors $v,u,w \in \MathR^k$,
 we say that $v$ is a {\em convex combination} of $u$ and $w$ if
there exists a real number $p$, $0 \le p \le 1$, such that $v =p\cdot
u + (1 - p) \cdot w$. 
A convex combination is called {\em trivial} if either $p = 0$ or $p =
 1$.
Otherwise, it is called {\em non-trivial}.
For a set $U \subseteq \MathR^k$ of vectors, we
say that $U$ is a {\em convexly independent set} if it contains no
three vectors $v, u,w \in U$ such that $v$ is a convex combination of
$u$ and $w$.
For a set $X \subseteq \MathR^k$ of vectors, 
the {\em exterior set} of
 $X$, 
denoted $\Ext(X)$, is the subset of $X$ that contains all vectors $v
\in X$ such that $v$ 
 cannot be expressed as a non-trivial convex combination of
vectors from $X$.

For a positive integer $\ell$, let $\beta_\ell$ denote the volume of
an $\ell$-dimensional ball of unit radius. It is well-known (see, e.g,
\cite{F82}, p.3)
that 
\begin{equation}
\label{eq:def_beta}
\beta_\ell ~=~ {{\pi^{\ell/2}} \over {\Gamma({\ell \over 2} + 1)}}~,
\end{equation}
where $\Gamma(\cdot)$ is the (Euler) Gamma-function.
We use the Gamma-function either with a positive integer parameter $n$ or with
a parameter $n + \half$ for a positive integer $n$.
In these cases the Gamma-function is given by
$\Gamma(n+1) = n!$ and 
\begin{equation}
\label{eq:Gamma}
\Gamma\left(n+ \half\right) ~=~ {{(2n)! \sqrt{\pi}} \over
  {2^{2n} n!}}~.
\end{equation}
(See \cite{F82}, p.178.)
 Observe also that
\begin{equation}
\label{eq:lb_Gamma}
\Gamma\left(n+ \half\right) ~=~ \left(n- \half\right)\left(n- {3 \over 2}\right) \cdot \ldots \cdot \half \cdot
\sqrt{\pi} ~\ge~ (n-1)! {{\sqrt{\pi}} \over 2}~.
\end{equation} 
By definition, it is easy to verify that for an integer $\ell$, $\ell
\ge 2$,
$\beta_\ell = \Theta({{\beta_{\ell - 1}} \over {\sqrt{\ell}}})$.

\section{Behrend Construction}
\label{sec:behrend}

The state-of-the-art lower bound for $\nu(n)$ due to Behrend \cite{B46}
states that for every positive integer $n$,
\begin{equation}
\label{eq:behrend_bnd}
\nu(n) ~=~ \Omega\left({n \over {{2^{2 \sqrt{2} \sqrt{\log n}}}
    \cdot \log^{1/4} n}}\right)~.
\end{equation}
In this paper we improve this bound by a factor of $\Theta(\sqrt{\log
  n})$,
and show that 
for every positive integer $n$,
\begin{equation}
\label{eq:our_bnd}
\nu(n) ~=~ \Omega\left({n \over {{2^{2 \sqrt{2} \sqrt{\log n}}}
    }} \cdot \log^{1/4} n\right)~.
\end{equation}
Note that it is sufficient to prove this bound only for all
sufficiently large values of $n$. The result for small values of $n$
follows by using a sufficiently small universal constant $c$ in the
definition of $\Omega$-notation.

We start with  a short overview of the original construction of
Behrend \cite{B46}.
Fix a sufficiently large positive integer $n$.
The construction involves a positive integer parameter $k$ that will
be determined later.
Set $y = n^{1/k}/2$.
In what follows we assume that $y$ is an integer. 
The case that $y$ is not an integer is analyzed later in the sequel.

Consider independent identically distributed random variables
$Y_1,Y_2,\ldots,Y_k$, with each $Y_i$ distributed uniformly over the
set $[\{0,y-1\}]$, for all $ i \in [\{k\}]$.
Set $Z_i = Y_i^2$, for all $i \in [\{k\}]$,
and $Z = \sum_{i=1}^k Z_i$.
It follows that for all  $i \in [\{k\}]$,
$$
\Expect(Z_i) ~ = ~ \sum_{j=0}^{y-1} {1 \over y} \cdot j^2 ~=~ {{y^2}
  \over 3} + \Theta(y)~.
$$
Let $\mu_Z = \Expect(Z)$ denote the expectation of the random variable
$Z$.
It follows that
\begin{equation}
\label{eq:muZ}
\mu_Z ~=~ {k \over 3} y^2 + \Theta(k \cdot y)~.
\end{equation}
Also, for all $i \in [\{k\}]$,
$\Var(Z_i) ~=~ \Expect(Z_i^2) - \Expect(Z_i)^2 ~=~ \Expect(Y_i^4) - 
{1 \over 9} y^4  + \Theta(y^3)$.
Hence
$$\Var(Z_i) ~=~ {{y^4} \over 5} + \Theta(y^3) - {{y^4} \over 9}  +
\Theta(y^3) 
~=~ {4 \over {45}} \cdot y^4 + O(y^3)~.$$
Hence 
$$\Var(Z) ~=~ k \cdot y^4 \cdot  {4 \over {45}} + O(k y^3) ~=~
 k \cdot y^4 \cdot  {4 \over {45}} \cdot ( 1 + O({1 \over y}))~,
$$ and the standard deviation of $Z$, $\sigma_Z$, satisfies
\begin{equation}
\label{eq:sigmaZ}
\sigma_Z ~=~ \sqrt{k} \cdot y^2 \cdot {2 \over {3 \cdot \sqrt{5}}}
\cdot (1 + O({1 \over y}))~.
\end{equation}
By Chebyshev inequality,
 for any $a > 0$,
$$\Prob(|Z - \mu_Z| > a \cdot \sigma_Z) \le {1 \over {a^2}}~.$$
Hence, for a fixed value of $a$, $a > 0$, at least $(1 - {1 \over
 {a^2}})$-fraction of all vectors $v$ from the set $[\{0,y-1\}]^k$
 have squared norm that satisfies
$$\mu_Z - a \cdot \sigma_Z \le ||v||^2 \le \mu_Z + a \cdot
 \sigma_Z~.$$
Note that each vector $v \in [\{0,y-1\}]^k$ has an integer squared
 norm.
By Pigeonhole Principle, there exists a value $T$ such that 
$\mu_z - a \cdot \sigma_Z \le T \le \mu_Z + a \cdot \sigma_Z$ that
 satisfies that at least $(1 - {1 \over {a^2}}) \cdot {1 \over {2a
 \cdot \sigma_Z}} \cdot y^k$ vectors from $[\{0,y-1\}]^k$ have squared
 norm $T$. Let $\cS$ denote the set of these vectors.
By (\ref{eq:sigmaZ}),
$$
|\cS| ~\ge ~ (1 - {1 \over {a^2}}) \cdot {1 \over {2a}}  {1 \over
  {\sqrt{k} \cdot y^2}} \cdot {{3 \sqrt{5}} \over 2} \cdot (1 - O({1
  \over y})) \cdot y^k ~=~ {{y^{k-2}}  \over {\sqrt{k}}} \cdot c~,$$
for a fixed positive constant $c = c(a)$.
Set $a = 2$. Now $c = c(2)$ is a universal constant, and
consequently,
$|\cS| = \Omega\left({{n^{k-2}} \over {2^k \sqrt{k}}}\right)$.
To maximize the right-hand-side, we set $k = \lceil \sqrt{2 \cdot \log
  n} \rceil$.
It follows that 
$$|\cS| ~=~ \Omega\left({n \over {{2^{2 \sqrt{2} \sqrt{\log n}}}
    \cdot \log^{1/4} n}}\right)~.$$
Observe that all vectors in $\cS$ have the same norm $\sqrt{T}$, and
thus, for every three vectors $v,u,w \in \cS$, $v \ne {{u + w} \over 2}$.
To obtain a progression-free set $S \subseteq [\{n\}]$ we consider
 coordinates of vectors from $\cS$ as digits of $(2y)$-ary
 representation.
Specifically, for every vector $v = (v_1,v_2,\ldots,v_k) \in \cS$, let 
$\hat{v} = \sum_{i=0}^{k-1} v_{i+1} \cdot (2y)^i$.
The set $S$ is now given by 
$S = \{\hat{v} \mid v \in \cS\}$.
Let $f(\cdot): \cS \rightarrow S$ denote this mapping.

Note that for every $v \in \cS$, 
$$0 ~<~ \hat{v} \le (2y)^k - 1 ~=~ n -1~.$$
Observe also that since 
$\cS \subseteq [\{0,y-1\}]^k$, 
the mapping $f$ is one-to-one, i.e., if $v \ne u$, $v,u
\in \cS$, then $\hat{v} \ne \hat{u}$.
Consequently, 
$$|S| ~=~ |\cS| ~=~ \Omega\left({n \over {{2^{2 \sqrt{2} \sqrt{\log n}}}
    \cdot \log^{1/4} n}}\right)~.$$
Finally, we argue that $S$ is a progression-free set. 
Suppose for contradiction that for three distinct numbers 
$\hat{v}, \hat{u}, \hat{w} \in S$,
$\hat{v} = {{\hat{u} + \hat{w}} \over 2}$. 
Let $u,v,w$ be the corresponding vectors in $\cS$, 
$v = (v_1,v_2,\ldots,v_k)$, $u = (u_1,u_2,\ldots,u_k)$,
$w = (w_1,w_2,\ldots,w_k)$.
Then
$$\hat{v} ~=~ \sum_{i=0}^{k-1} {{u_{i+1} + w_{i+1}} \over 2} \cdot
(2y)^i ~=~ \sum_{i=0}^{k-1} v_{i+1} \cdot (2y)^i~.$$
However, since all the coordinates
$v_1,v_2,\ldots,v_k,u_1,u_2,\ldots,u_k, w_1,w_2,\ldots,w_k$ are in
$[\{0,y-1\}]$, it follows that $v_i = {{u_i + w_i} \over 2}$, for
every index $i \in [\{k\}]$.
Consequently, $v = {{u+w} \over 2}$, a contradiction to the assumption
that $||v|| = ||u|| = ||w||$.
Hence $S$ is a progression-free set of size 
$\Omega({n \over {{2^{2 \sqrt{2} \sqrt{\log n}}}
    \cdot \log^{1/4} n}})$.

Consider now the case that $y = {{n^{1/k}} \over 2}$ is not an integer
number.
In this case the same construction is built with $\lfloor y \rfloor$
instead of $y$.
Set $n' = (2 \lfloor y \rfloor)^k$. 
By previous argument, we obtain a progression-free set $S$ that satisfies 
$$|S| ~=~ \Omega\left({n' \over {{2^{2 \sqrt{2} \sqrt{\log n'}}}
    \cdot \log^{1/4} n'}}\right) ~=~
\Omega\left({n' \over {{2^{2 \sqrt{2} \sqrt{\log n}}}
    \cdot \log^{1/4} n}}\right)~.$$
Observe that 
${n \over {n'}} ~\le~ \left({y \over {y-1}}\right)^k ~=~ 1 +
\Theta({k \over y}) ~=~ 1 + \Theta\left({{\sqrt{\log n}} \over
  {2^{(1/\sqrt{2}) \cdot \sqrt{\log n}}} }\right)$.

Hence 
$ |S| ~=~ \Omega\left({n \over {{2^{2 \sqrt{2} \sqrt{\log n}}}
    \cdot \log^{1/4} n}}\right)$,
and we are done.

\section{Our Construction}
\label{sec:our}

In this section we present our construction of progression-free sets
$S \subseteq [\{n\}]$
 with at least $\Omega\left({n \over {{2^{2 \sqrt{2} \sqrt{\log n}}}
    }} \cdot \log^{1/4} n\right)$ elements.
Fix $k = \lceil \sqrt{2 \log n} \rceil$, and $y = n^{1/k}/2$.
Observe that 
\begin{equation}
\label{eq:y_def}
{{2^{k/2}} \over {2 \sqrt{2}}} ~=~ {1 \over {2 \sqrt {2}}} \cdot 
2^{{\sqrt{\log n} \over \sqrt{2}}} ~\le~ y ~\le~ \half \cdot  2^{{\sqrt{\log n} \over \sqrt{2}}} ~=~ {{2^{k/2}}
  \over 2}~.
\end{equation}
For convenience 
we assume that $y$ is an integer. If this is not the
case, the same analysis applies with  minor adjustments. 
(Specifically, we set $y = \lfloor  n^{1/k}/2 \rfloor$. By the same
argument as we used in Section \ref{sec:behrend}, the resulting lower
bound will be at most by a constant factor smaller than in the case when
$n^{1/k}/2$ is an integer.)

Consider the $k$-dimensional ball  centered at the origin that has
radius $R'$ given by 
\begin{equation}
\label{eq:Rdef}
R'^2 ~=~ \mu_Z ~=~ {k \over 3} y^2 + \Theta(k y)~.
\end{equation}
(See (\ref{eq:muZ}).) By Chebyshev inequality, the  annulus
$\hcS$ of all vectors with squared norm in $[R'^2 - 2 \cdot \sigma_Z,
  R'^2  + 2 \cdot \sigma_Z]$ contains at least ${3 \over 4} \cdot y^k$
integer points of the discrete cube $C = [\{0,y-1\}]^k$.

Fix a parameter $g = \eps \cdot k$, for a universal constant $\eps >
0$ that will be determined later.
Partition the annulus $\hcS$ into $\lceil {{4 \sigma_Z} \over g} \rceil
= \ell$ annuli $\hcS_1,\hcS_2,\ldots,\hcS_\ell$, with the annulus
$\hcS_i$ containing all vectors with squared norms in the range 
$[R'^2 - 2  \sigma_Z + (i-1) \cdot g, R'^2 - 2  \sigma_Z + i
  \cdot g)$, for $i \in [\{\ell - 1\}]$, and 
$[R'^2 - 2 \sigma_Z + (\ell - 1) \sigma_Z, R'^2 + 2 \sigma_Z]$ for $i
  = \ell$.

Observe that for distinct indices $i,j \in [\{\ell\}]$, the sets of
integer points in $\hcS_i$ and $\hcS_j$ are disjoint. Thus,
by the Pigeonhole Principle, there exists an index $i \in [\{\ell\}]$
such that the annulus $\hcS_i$ contains at least 
\begin{equation}
\label{eq:num_pnts}
{3 \over {4 \ell}}
  \cdot y^k ~=~ \Omega(g \cdot {{y^{k-2}} \over {\sqrt{k}}}) ~=~
  \Omega(\eps \sqrt{k} \cdot y^{k-2})
\end{equation}
 integer points of $C \cap \hcS$.
In other words, there exists a radius $R$, 
$R^2 \in [R'^2 - 2\sigma_Z,R'^2 + 2 \sigma_Z]$, such that the
 annulus $\cS$ that contains all vectors with squared norm in
the range $[R^2 -g,R^2]$ contains at least $\Omega(\sqrt{k} \cdot
y^{k-2})$ integer points of $C \cap \hcS$.

By (\ref{eq:muZ}), (\ref{eq:sigmaZ}), and (\ref{eq:Rdef}),
\begin{equation}
\label{eq:ubR}
R^2 ~\le~ R'^2 + 2 \sigma_Z ~\le~ {k \over 3} \cdot y^2 + O(k \cdot y)
+ O(\sqrt{k} \cdot y^2) ~\le ~ {k \over 3} \cdot y^2 \left(1 +
O\left({1 \over {\sqrt{k}}}\right)\right)~.
\end{equation}

Let $\tcS$ be the set of integer points of $C \cap \cS$.
We will show that that $\tcS$ contains a convexly independent subset 
$\ccS$ with at least $|\ccS| \ge {{|\tcS|} \over 2}$ integer points.
Consequently, 
\begin{equation}
\label{eq:ccS}
|\ccS| ~\ge~ 
{{|\tcS|} \over 2} ~=~ \Omega(\sqrt{k} \cdot y^{k-2}) ~=~
\Omega\left(\log^{1/4} n \cdot {n \over {2^{2\sqrt{2} \sqrt{\log n}}}}\right)~.
\end{equation}
Consider the set $\chS = f(\ccS)$ constructed from $\ccS$ by the
mapping $f$ described in Section \ref{sec:behrend}.
Since $\cS$ is a convexly independent set, by the same argument as in
Section \ref{sec:behrend}, $|\chS| = |\ccS|$, and moreover, 
$\chS$ is a progression-free set. Hence 
$|\chS| = \Omega\left(\log^{1/4} n \cdot {n \over {2^{2\sqrt{2}
      \sqrt{\log n}}}}\right)$, and our result follows.

The following lemma is useful for showing an upper bound on the number
of integer points in $\cS$ that do not belong to the exterior set of
$\tcS$,
$\Ext(\tcS)$. This lemma is due to Coppersmith \cite{C03}.

Let $\cB = \cB(R,0)$ denote the $k$-dimensional ball of radius $R$
centered at the origin, and $B = B(R,0)$ denote the set of integer
points contained in this ball. Denote $T = R^2$.

\begin{lemma}
\label{lm:wipe} \cite{C03}
Let  $b \in B \setmns \Ext(B)$ be an integer point that satisfies
$T - g \le ||b||^2 \le T$. Then there exists a non-zero integer vector
$\delta$  that satisfies $ 0 \le \langle b, \delta \rangle \le g$ and
$0 < || \delta||^2 \le g$.
\end{lemma}
\proof
Since $b \in B \setmns \Ext(B)$, there exist two integer points $a$
and $c$ in $B$ and a constant $p$, $0 < p < 1$, such that
$b = p \cdot a + (1 - p) \cdot c$.
Since $a,c \in B$, $||a||^2, ||c||^2  \le T$.
Observe that either $\langle a, b \rangle$ or $\langle c,b \rangle$
is greater or equal than $||b||^2$. (Otherwise,
$||b||^2 = \langle p a + (1 -p)c, b \rangle = p \cdot \langle a,b
\rangle + (1- p) \cdot \langle c,b \rangle < ||b||^2$,
contradiction.)

Suppose without loss of generality that $\langle a, b \rangle \ge
||b||^2$.
Then $\langle a - b,b \rangle \ge 0$. Set $\delta = a - b$.
Since $a,b \in B$ are integer points, it follows that $\delta$ is an
integer point as well. Moreover, since $0 < p < 1$, we have $\delta \ne 0$.
Moreover, 
$$T ~\ge~ ||a||^2 ~=~ ||b + \delta||^2 ~=~ ||b||^2 + 2 \langle
b,\delta \rangle + ||\delta||^2~.$$
Recall that $||b||^2 \ge T - g$.
Hence $2\langle
b,\delta \rangle + ||\delta||^2 \le g$.
As $\langle b,\delta \rangle = \langle a - b,b \rangle \ge 0$, it
follows that $\langle b,\delta \rangle, ||\delta||^2 \le g$, as required.
\QED

Observe that $\delta \in \MathR^k$ is an integer vector, and
$||\delta||^2 \le g = \eps \cdot k$. Consequently, the vector $\delta$
may contain at most $\eps \cdot k$ non-zero entries.
This property will be helpful for our argument.

Denote the number of integer vectors $\delta$ that have squared norm 
at most $g$ by $\hD(g)$. The next lemma provides an upper bound on
$\hD(g)$.

\begin{lemma}
\label{lm:few}
For any $\eps > 0$ there exists $\eta = \eta(\eps) > 0$ such that 
$\lim_{\eps \to 0} \eta(\eps) = 0$, and $\hD(g) = O(2^{\eta \cdot k})$.
\end{lemma}
\proof
Fix an integer value $h$, $1 \le h \le g$. 
First, we count the number $N(h)$ of $k$-tuples $(q_1,q_2,\ldots,q_k)$
of non-negative integer numbers that sum up to $h$.

Consider permutations of $(k - 1 + h)$ elements of two types, with $h$
elements of the first type and $k-1$ elements of the second type.  Elements
of the first type are called ``balls'', and elements of the second
type are called ``boundaries''. Two permutations $\sigma$ and
$\sigma'$ are said to be {\em equivalent} if they can be obtained one
from another by permuting balls among themself, and permuting
boundaries among themself. 

Let $\Pi$ be the induced equivalence relation.
Observe that there is a one-to-one mapping
between $k$-tuples $(q_1,q_2,\ldots,q_k)$ of non-negative 
integer numbers that sum
up to $h$ and the equivalence classes of the relation $\Pi$.
Hence $N(h)$ is equal to the number of equivalence classes of $\Pi$,
i.e.,
$$N(h) ~=~ {{(k - 1 + h)!} \over {(k-1)! \cdot h!}} ~=~ {{k - 1 + h} \choose
  h}~.$$
In a $k$-tuple $(\delta_1,\delta_2,\ldots,\delta_k)$ of integer
numbers such that $\sum_{i=1}^k \delta_i^2 = h$, there can be at most
$h$ non-zero entries. Hence, for a fixed $k$-tuple of integers
$(q_1,q_2,\ldots,q_k)$ such that $\sum_{i=1}^k q_i = h$, there may be
at most $2^h$ $k$-tuples $(\delta_1,\delta_2,\ldots,\delta_k)$ of
integers such that $\delta_i^2 = q_i$ for every index $i \in [\{k\}]$.
Thus, the overall number $D(h)$ of integer $k$-tuples
$(\delta_1,\delta_2,\ldots,\delta_k)$ such that $\sum_{i=1}^k
\delta_i^2 = h$ satisfies 
$$D(h) ~\le~ 2^h \cdot N(h)  = 2^h {{k - 1 + h} \choose
  h}~.$$
Note that
${{k - 1 + h} \choose
  h} \le {{k - 1  + g} \choose g}$, for every integer $h$, $1 \le h
\le g$.  Hence the number $\hD(g)$ of integer $k$-tuples
$(\delta_1,\delta_2,\ldots,\delta_k)$ with $1 \le \sum_{i=1}^k
\delta_i^2 \le g$ satisfies
\begin{eqnarray*}
\hD(g) &=& \sum_{h=1}^g D(h) ~\le~ \sum_{h=1}^g 2^h \cdot N(h) ~\le~
N(g) \cdot 2^{g+1} ~\le~ 2^{g+1} \cdot {{k  + g } \choose g} \\
& \le &
2^{g+1} \left({{e(k+g)} \over g}\right)^g ~=~
2 \cdot (2 e)^g \left(1 + {1 \over \eps}\right)^{\eps \cdot k} ~=~
2 \cdot 2^{(\log 2e + \log (1 + {1 \over \eps})) \eps \cdot k}~.
\end{eqnarray*}
Denote $\eta = \eta(\eps) = \eps (\log 2e + \log (1 + {1 \over
  \eps}))$.
Then $\hD(g) \le 2 \cdot 2^{\eta(\eps) \cdot k}$.
Finally, 
$$\lim_{\eps \to 0} \eta(\eps)~=~ 
\lim_{\eps \to 0} {{\log(1 + {1 \over \eps})} \over {1 \over \eps}} ~=~
{1 \over {\ln 2}} \cdot \lim_{y \to \infty}
 {{\ln (1 + y)}
  \over y} ~=~ 0~,$$ completing the proof.
\QED

Consider again the  annulus $\cS = \{\alpha \in \MathR^k \mid T
- g \le ||\alpha||^2 \le T\}$, and the set $\tcS$ of integer points of
$\cS$.
For an integer vector $\delta$ that satisfies $0 < ||\delta||^2 \le
g$, let $\hZ(\delta)$ denote the set of integer points $ b \in \tcS$
that satisfy $0 \le \langle b,\delta \rangle  \le g$.
Let $\hW(\delta) = \hZ(\delta) \cap C$ denote the intersection of
$\hZ(\delta)$ with the discrete cube $C = [\{0,y-1\}]^k$, and let $W(\delta) =
|\hW(\delta)|$. Also, let $\hW = \bigcup \{\hW(\delta) \mid 0 <
||\delta||^2 \le g\}$, and $W = |\hW|$.

Let $\hN$ denote the set of integer points of $C \cap \cS$ that do not
belong to $\Ext(B)$, and $N = |\hN|$.
By Lemma \ref{lm:wipe},
$\hN \subseteq \hW$, and consequently,
\begin{equation}
\label{eq:N_bnd}
N ~\le ~ W ~\le ~ \sum \{ W(\delta) \mid 0 < ||\delta||^2 \le g \}~.
\end{equation}
Fix a vector $\delta$, $0 < ||\delta||^2 \le g$.
In the sequel we provide an upper bound for $W(\delta)$.

Observe that since $\hW(\delta)$ is a set of integer points, it
follows that for every $b \in \hW(\delta)$,
$\langle b,\delta \rangle \in [\{0,g\}]$.

For an integer number $h \in [\{0,g\}]$, let $\hW(\delta,h)$ denote
the subset of $\hW(\delta)$ of  integer points $b$ that satisfy
$\langle b,\delta \rangle = h$. Let $W(\delta,h) = |\hW(\delta,h)|$.
Observe that for distinct values $h \ne h'$, $h,h' \in [\{0,g\}]$, the
sets $\hW(\delta,h)$ and $\hW(\delta,h')$ are disjoint.
Consequently,
\begin{equation}
\label{eq:W}
W(\delta) ~ = ~ \sum_{h= 0}^g W(\delta,h)~.
\end{equation}
Next, we provide an upper bound for $W(\delta,h)$.

Consider the hyperplane $\cH = \{\alpha \in \MathR^k \mid \langle
\alpha,\delta \rangle = h\}$.
Observe that $\hW(\delta,h) = \cH \cap \cS \cap C$ is the intersection
of the hyperplane $\cH$ with the  annulus $\cS$ and with the
discrete cube C.

Let $S$ denote the $k$-dimensional sphere with squared radius $T$
centered at the origin, i.e., $S = \{ \alpha \in \MathR^k \mid
||\alpha||^2 = T\}$. Consider the intersection $S'$ of $S$ with the
hyperplane $\cH$.

\begin{lemma}
\label{lm:intersect}
$S' \subseteq \cH$ is a $(k-1)$-dimensional sphere with squared radius 
$(T - {{h^2} \over {||\delta||^2}})$ centered at 
${ h \over {||\delta||^2}} \cdot \delta$.
\end{lemma}
\proof
For a vector $\alpha \in S \cap \cH$, 
$$
||\alpha - { h \over {||\delta||^2}} \cdot \delta ||^2 ~ = ~
\sum_{i=1}^k (\alpha_i - {h \over {||\delta||^2}} \cdot \delta_i)^2
~=~
||\alpha||^2 + {{h^2} \over {||\delta||^2}} - 2 {h \over
  {||\delta||^2}} \langle \alpha, \delta \rangle ~=~
 ||\alpha||^2 - {{h^2} \over {||\delta||^2}}~.
$$
(For the last equality, note that since $\alpha \in \cH$, we have
 $\langle \alpha, \delta
\rangle = h$.)
\QED

Recall that for a vector $\alpha \in \cS$, $T  -g \le ||\alpha||^2 \le
T$. Hence
 the intersection of the hyperplane $\cH$ with the
 annulus $\cS$ is the $(k-1)$-dimensional  annulus $\cS'
\subseteq \cH$, centered at ${ h \over {||\delta||^2}} \cdot \delta$,
containing vectors $\alpha$ such that 
$$T - g - {{h^2} \over {||\delta||^2}} ~\le~ ||\alpha - { h \over
  {||\delta||^2}} \cdot \delta ||^2 ~\le~ T - {{h^2} \over
  {||\delta||^2}}~.$$
Let $T ' = T - {{h^2} \over {||\delta||^2}}$.
Then $\cS'$ is given by
$$\cS' ~=~ \{\alpha \in \cH \mid T' - g \le ||\alpha - { h \over
  {||\delta||^2}} \cdot \delta ||^2  \le T'\}~.$$
Note that since $h \ge 0$, $T' \le T$ for all $h$ and $\delta$.

Recall that our goal at this stage is to provide an upper bound for
the number $W(\delta,h)$ of integer points in $\hW(\delta,h) = \cH
\cap \cS \cap C = \cS' \cap C$.
Let $\cC = [0,y-1]^k$ be the (continuous) cube. (The discrete cube 
$C = [\{0,y-1\}]^k$ is the set of integer points of $\cC$.)
Let $\tW = \cS' \cap \cC$.
Since $\hW(\delta,h)$ is the set of integer points in $\tW$, we are
interested in providing an upper bound for the number of integer
points in $\tW$. Our strategy is to show an upper bound for the
$(k-1)$-dimensional volume $\Vl(\tW)$ of $\tW$, and to use standard
estimates for the discrepancy between $\Vl(\tW)$ and the number of
integer points in $\tW$.

Let $\cH' = \{\alpha \in \MathR^k \mid \langle \alpha,\delta \rangle =
0\}$ be the parallel hyperplane to $\cH$ that passes through the
origin.
Next, we construct an orthonormal basis $\Upsilon =
\{\gamma_1,\gamma_2,\ldots,\gamma_{k-1}\}$ for $\cH'$. This basis will
be useful for estimating $\Vl(\tW)$.

Recall that $\delta$ satisfies $0 < ||\delta||^2 \le g = \eps \cdot
k$, and it is an integer vector. Consequently, $\delta =
(\delta_1,\delta_2,\ldots,\delta_k)$ contains at most $g = \eps \cdot
k$ non-zero entries. Let $I \subseteq [\{k\}]$ be the subset of
indices such that $\delta_i \ne 0$. Let $m = |I|$. 
It follows that $m  \le g = \eps \cdot
k$.

For every vector $\alpha = (a_1,a_2,\ldots,a_k) \in \cH'$, it holds
that
\begin{equation}
\label{eq:orthog}
\sum_{i \in I} a_i \delta_i = 0~.
\end{equation}
Let $\gamma^{(1)}, \gamma^{(2)},\ldots,\gamma^{(m-1)}$ be an arbitrary
orthonormal basis for the solution space of the equation
(\ref{eq:orthog}).
These vectors are in $\MathR^m$.
For each index $j \in [\{m-1\}]$, we view the vector $\gamma^{(j)}$ as 
$\gamma^{(j)} = (\gamma_i^{(j)} \mid i \in I)$.

We form orthonormal vectors
$\hgamma^{(1)},\hgamma^{(2)},\ldots,\hgamma^{(m-1)} \in \MathR^k$ in
the following way.
For each index $j \in [\{m-1\}]$, and each index $i \in I$, the $i$th
entry $\hgamma_i^{(j)}$ of $\hgamma^{(j)}$ is set as $\gamma_i^{(j)}$,
and for each index $i \in [\{k\}] \setmns I$, the entry
$\hgamma_i^{(j)}$ is set as zero.
Also, for each index $i \in [\{k\}] \setmns I$, we insert the vector 
$\xi_i = (0,0,\ldots,0,1,0,\ldots,0)$, $\xi_i \in \MathR^k$, with 1 at
the $i$th entry and zeros in all other entries into the basis
$\Upsilon$.
Observe that $\xi_i \in \cH'$.
The resulting basis $\Upsilon$ is
$\{\hgamma^{(1)},\hgamma^{(2)},\ldots,\hgamma^{(m-1)}\} \cup
\{ \xi_i \mid i \in [\{k\}] \setmns I\}$.
It is easy to verify that $\Upsilon$ is an orthonormal basis for $\cH'$.

Order the vectors of $\Upsilon$ so that $\hgamma^{(j)} = \gamma_j$ for
all
$j \in [\{0,m-1\}]$, and 
so that the  vectors $\{\xi_i \mid i \in [\{k\}] \setmns
I\}$ appear in an arbitrary order among
$\gamma_m,\gamma_{m+1},\ldots,\gamma_{k-1}$.

Move the origin to the center ${h \over {||\delta||^2}} \cdot \delta$
of the  annulus $\cS'$, and rotate the annulus so that new axes
become the colinear with vectors
$\gamma_1,\gamma_2,\ldots,\gamma_{k-1}$ of the orthonormal basis
$\Upsilon$. Obviously, this mapping is volume-preserving.

For a vector $\zeta \in \cH'$, let
$\zeta_1[\Upsilon],\zeta_2[\Upsilon],\ldots,\zeta_{k-1}[\Upsilon]$
denote the coordinates of $\zeta$ with respect to the basis $\Upsilon$,
i.e., $\zeta_i[\Upsilon] = \langle \zeta - {h \over {||\delta||^2}} \cdot
\delta, \gamma_i \rangle$.
Observe that since $\langle \delta,\gamma_i \rangle = 0$ for all $i \in
[\{k-1\}]$,
it follows that $\zeta_i[\Upsilon] =  \langle \zeta,\gamma_i \rangle$,
for all $i \in
[\{k-1\}]$.

\begin{lemma}
\label{lm:many_nonneg}
For  a vector $\zeta \in \tW = \cS' \cap \cC$,
and an index $i \in [\{m,\ldots,k-1\}]$, we have
$\zeta_i[\Upsilon] \ge 0$. In particular,
$\zeta$ has at least $(1-\eps)\cdot k$ non-negative
coordinates with respect to the basis $\Upsilon$.
\end{lemma}
\proof
Note that for every index $i \in [\{m,k-1\}]$, all entries of
$\gamma_i$ are non-negative. (Because these are the vectors $\xi_j$,
$j \in  [\{k\}] \setmns
I\}$ of the standard Kronecker basis.)
Since $\zeta \in \cC = [0,y-1]^k$, it follows that for all indices
$i \in [\{m,k-1\}]$, the $i$th coordinate of $\zeta$ with respect to
the basis $\Upsilon$ is non-negative,
that is, $\zeta_i[\Upsilon] = \langle \zeta,\gamma_i \rangle \ge 0$.
Hence $\zeta$ has at least $(k-1) -(m-1) \ge (1-\eps) \cdot k$ non-negative
coordinates with respect to the basis $\Upsilon$.
\QED

Recall that $\tW \subseteq \cS'$, and the  annulus $\cS'$ is
given by (with respect to the basis $\Upsilon$)
$$\cS' ~=~ \{ \alpha \in \MathR^{k-1} \mid T' -g \le || \alpha||^2 \le
T'\}~.$$
Let 
$$\tQ ~=~ \{\alpha \in \MathR^{k-1} \mid T' - g \le ||\alpha||^2 \le T',
~\forall i \in [\{m,k-1\}], \alpha_i[\Upsilon] \ge 0\}$$
be the intersection of $\cS'$ with the $(k - m)$ half-spaces
$\alpha_i[\Upsilon] \ge 0$, for all $i \in [\{m,k-1\}]$. 
Let $\cS''$ be the intersection of the  annulus $\cS'$ with the
positive octant (with respect to $\Upsilon$), i.e.,
$$\cS'' = \{\alpha \in (\MathR^+)^{k-1} \mid T' - g \le ||\alpha||^2
\le T'\}~.$$
It follows that $\tW \subseteq \tQ$, and 
$\Vl(\tW) \le \Vl(\tQ) = 2^m \cdot \Vl(\cS'') ~\le 
2^{\eps \cdot k - 1} \cdot
\Vl(\cS'')$.

Next, we provide an upper bound for $\Vl(\cS'')$.

\begin{lemma}
\label{lm:vol}
For a sufficiently large integer $k$,
$$\Vl(\cS'') ~\le~ g \cdot \left({{\pi e} \over 6}\right)^{k/2} \cdot
y^{k-3} \cdot 2^{O(\sqrt{k})}~.$$
\end{lemma}
\proof
Let $R' = \sqrt{T'}$. Observe that $R' \le R = \sqrt{T}$.
Let $\beta_{k-1}$ be the volume of the $(k-1)$-dimensional ball of
unit radius.
Then
$\Vl(\cS'') = {1 \over {2^{k-1}}} \beta_{k-1} ((T')^{{k-1} \over 2} - 
(T' - g)^{{k-1} \over 2})$.
Note that
\begin{equation}
 (R'^2 - g)^{{k-1} \over 2} ~ = ~ \left(1 - {g \over
    {R'^2}}\right)^{{k-1} \over 2} \cdot R'^{k-1} 
~\ge ~R'^{k-1} \left(1 - {{g(k-1)} \over {2 R'^2}}\right) 
~ \ge ~ R'^{k-1} - R'^{k-3} g \cdot k~.
\end{equation}
Hence by (\ref{eq:def_beta}), 
\begin{equation}
\label{eq:volS}
\Vl(\cS'') ~\le~ {1 \over {2^{k-1}}} \cdot k \cdot g \cdot \beta_{k-1} \cdot
R'^{k-3} ~ \le ~  {1 \over {2^{k-1}}} \cdot k \cdot g \cdot {{\pi^{{k-1} \over
      2}} \over {\Gamma\left({{k+1} \over 2}\right)}} \cdot R^{k-3}~.
\end{equation}
By (\ref{eq:ubR}),
 $T = R^2 \le {k \over 3} \cdot y^2 \left(1 + O\left({1
  \over {\sqrt{k}}}\right)\right)$, and so
$R \le \sqrt{k \over 3} \cdot y 
\left(1 + O\left({1
  \over {\sqrt{k}}}\right)\right)$.
Hence
$$\Vl(\cS'') ~\le~ {1 \over {2^{k-1}}} \cdot k \cdot g \cdot {{\pi^{{k-1} \over
      2}} \over {\Gamma\left({{k+1} \over 2}\right)}} \cdot 
\left({k \over 3}\right)^{{k-3} \over 2} \cdot y^{k-3} \cdot
      2^{O(\sqrt{k})}~.$$
By Stirling formula, if $k+1$ is even then
for a sufficiently large $k$,
\begin{eqnarray*}
\Gamma\left({{k+1} \over 2}\right) &\ge & \left({{k-1} \over 2}
\right)! ~\ge~ \sqrt{{k-1} \over 2} \cdot {{\left({{k-1} \over
      2}\right)^{{k-1} \over 2}} \over {e^{{k-1} \over 2}}} \\
& = & e^{1/2} {{{{\left({k \over 2}\right)^{k \over 2} \left(1 - {1
      \over k}\right)^{k \over 2}}}} \over {e^{k/2}}}~\ge~
\half \cdot {{k^{k \over 2}} \over {(2e)^{k \over 2}}}~.
\end{eqnarray*}
By (\ref{eq:lb_Gamma}),
 if $k+1$ is odd then for a sufficiently large $k$,
\begin{eqnarray*}
\Gamma\left({{k+1} \over 2}\right) & = &  \Gamma\left({k \over 2} +
\half\right)
~\ge~ {{\sqrt{\pi}} \over 2} \left({k \over 2} -1\right)!  ~\ge~
{{\pi e} \over {\sqrt{2}}} \cdot \sqrt{{k \over 2} - 1} \cdot
{{\left({k \over 2} - 1\right)^{{k \over 2} - 1}} \over 
{e^{{k \over 2}  - 1}}} \\
& \ge &
 {{\pi e} \over {\sqrt{k}}} \cdot {{\left({k \over 2}\right)^{k \over
       2} \cdot \left(1 - {2 \over k}\right)^{k \over 2}} \over {e^{k
       \over 2}}} ~\ge ~
{{1} \over {\sqrt{k}}} \cdot {{k^{k/2} } \over {(2e)^{k/2}}}~.
\end{eqnarray*}
Hence in both cases, for a sufficiently large $k$,
$$\Gamma\left({{k + 1} \over 2}\right) ~\ge ~ 
{1 \over {\sqrt{k}}} \cdot {{k^{k \over 2}} \over {(2e)^{k \over 2}}}~.$$
Consequently, 
\begin{eqnarray*}
\Vl(\cS'') &\le & (k \cdot g) \cdot {1 \over {2^{k-1}}} \cdot
{{\pi^{k \over 2} \sqrt{k} \cdot (2e)^{k \over 2}} \over {\sqrt{\pi} \cdot k^{k
      \over 2}}} \cdot {{k^{{k-3} \over 2}} \over {3^{{k-3} \over 2}}}
\cdot y^{k-3} \cdot 2^{O(\sqrt{k})}  \\
& = &
O(1) \cdot (k \cdot g) \cdot k^{{-3} \over 2} \sqrt{k} \cdot
\left({{\pi e} \over 6}\right)^{k \over 2} \cdot y^{k-3} \cdot
2^{O(\sqrt{k})} 
~\le ~ g \cdot \left({{\pi e} \over 6}\right)^{k \over 2} \cdot y^{k-3} \cdot
2^{O(\sqrt{k})}~.
\inQED
\end{eqnarray*}

We conclude that
\begin{equation}
\label{eq:vlW}
\Vl(\tW) ~\le~ \Vl(\tQ) ~\le~
2^{\eps k -1} \cdot \Vl(\cS'') ~\le~
\half \cdot g \cdot 2^{\eps k} \cdot \left({{\pi e} \over 6}\right)^{k \over 2} \cdot y^{k-3} \cdot
2^{O(\sqrt{k})}~.  
\end{equation}

Since $\tW \subseteq \tQ$, the number $W(\delta,h)$ of integer points
in $\tW$ is at most the number $Q$ of integer points in $\tQ$.
In Section \ref{sec:int_pnt} we will show that $Q$ is not much larger
than $\Vl(\tQ)$. Specifically, 
\begin{equation}
\label{eq:Q_bnd}
Q ~\le ~ k^{O(1)} \cdot 2^{\eps k} \cdot \left({{\pi e} \over 6}\right)^{k \over 2} \cdot y^{k-3} \cdot
2^{O(\sqrt{k})}~.
\end{equation}
We remark that this estimate is quite crude, as it says that the
number $Q$ of integer points in $\tQ$ cannot be larger than by a
factor of $k^{O(1)}$ than $\Vl(\tQ)$. However, it is sufficient for
our argument. 

By (\ref{eq:Q_bnd}),
\begin{equation}
\label{eq:W_bnd}
W(\delta,h) ~\le~ Q ~\le~ k^{O(1)} \cdot 2^{\eps k} \cdot \left({{\pi e} \over 6}\right)^{k \over 2} \cdot y^{k-3} \cdot
2^{O(\sqrt{k})}~.
\end{equation}
By (\ref{eq:W}), 
$$W(\delta) ~ = ~ \sum_{h= 0}^g W(\delta,h) ~\le~
(g + 1) \cdot k^{O(1)} \cdot 2^{\eps k} \cdot \left({{\pi e} \over 6}\right)^{k \over 2} \cdot y^{k-3} \cdot
2^{O(\sqrt{k})}~.$$

Hence by (\ref{eq:N_bnd}), the overall number $N$ of integer points in
$C \cap \cS$ that do not belong to $\Ext(B)$ (and thus, do not belong
to $\Ext(C \cap \cS)$, because $\cS \subseteq B$) satisfies
$$N ~\le ~ \sum_{0 < ||\delta||^2 \le g} W(\delta) ~\le~
k^{O(1)} \cdot 2^{\eps k} \cdot \left({{\pi e} \over 6}\right)^{k
  \over 2}
 \cdot y^{k-3} \cdot
2^{O(\sqrt{k})} \cdot \hD(g)~.$$
Recall that $g \le k$.
By Lemma \ref{lm:few}, and since for a sufficiently
large $k$, 
$k^{O(1)} \le 2^{O(\sqrt{k})}$, it follows that
$$
N ~\le ~k^{O(1)} \cdot 2^{\eps k} \cdot \left({{\pi e} \over 6}
\right)^{k \over 2} \cdot y^{k-3} \cdot
2^{O(\sqrt{k})} \cdot O(2^{\eta \cdot k})
~=~ 2^{((\eps + \eta(\eps) + O(1/ \sqrt{k})) + \half \log {{\pi e}
    \over 6}) \cdot k} \cdot y^{k-3}~.
$$
By (\ref{eq:num_pnts}), the set $\tcS$ of integer points of $C \cap \cS$
contains
$$|\tcS| ~ = ~ \Omega(\eps \sqrt{k} \cdot y^{k-2})$$
integer points.
By (\ref{eq:y_def}),
\begin{equation}
\label{eq:y_large}
y ~=~ {{2^{k/2}} \over 2} ~>~  
2 \cdot 
O\left({1 \over {\eps \cdot \sqrt{k}}}\right)
2^{((\eps + \eta(\eps) + O(1/ \sqrt{k})) + \half \log {{\pi e}
    \over 6}) \cdot k}
\end{equation}
whenever $\eps > 0$ and $k$ satisfy
$$1 ~>~ \log {{\pi e} \over 6} + (\eps + \eta(\eps)) + O\left({1 \over
  {\sqrt{k}}}\right)~.$$
By Lemma \ref{lm:few}, $\lim_{\eps \to 0} \eta(\eps) = 0$.
Thus, for a sufficiently small universal constant $\eps > 0$, and
  sufficiently large $k$, 
the inequality (\ref{eq:y_large}) holds, and thus 
$|\tcS| ~\ge~ 2 N$.
(More specifically, one needs to set $\eps$ so that
$0 < \eps + \eta(\eps) < 1 - \log {{\pi e} \over 6}$.)
Hence the set $\tcS$ contains a subset $\ccS$ of integer points that
  belong to $\Ext(B)$, and moreover,
$$|\ccS| ~\ge~ |\tcS| - N ~\ge ~ \half |\tcS| ~=~
\Omega(\eps \cdot \sqrt{k} \cdot  y^{k-2}) ~=~
\Omega(\log^{1/4} n \cdot {n \over {2^{2\sqrt{2} \sqrt{\log n}}}})~.$$

\section{ Discrepancy between Volume and \\ Number of Integer Points}
\label{sec:int_pnt}

Consider the  annulus $\cS' = \{ \alpha =
(\alpha_1,\alpha_2,\ldots,\alpha_{k-1}) \in \MathR^{k-1} \mid T' - g
\le ||\alpha||^2 \le T'\}$, and its intersection $\tQ$ with the
half-spaces
$\alpha_i \ge 0$ for all $i \in [\{m,k-1\}]$.
In this section we argue that the number $Q$ of integer points in
$\tQ$ is not much larger than $\Vl(\tQ)$. 
Specifically,
we show that 
\begin{equation}
\label{eq:Q}
Q ~=~ 2^{O(\sqrt{k})} \cdot 2^{\eps k} \cdot \left({{\pi e} \over
  6}\right)^{k/2} \cdot y^{k-3}~.
\end{equation}
This proves (\ref{eq:Q_bnd}), and hence completes the proof of our lower
bound.

Consider the $(k-1)$-dimensional ball $B$ of squared radius $t$
centered at the origin, for some sufficiently large $t > 0$.
Let $A(B)$ denote the number of integer points in $B$.
For a positive
integer $j$, let $V_j(t)$ denote the volume of the $j$-dimensional
ball of squared radius $t$ centered at the origin.
It is well-known (see, e.g., the survey of Adhikari \cite{Adh94}) that for
a constant dimension $k$,
$|A(B) - V(B)| = O(V_{k-3}(t))$.
However, in our case the dimension $k$ {\em grows logarithmically} with $t$.
Fortunately, the following analogous inequality holds in this case:
\begin{equation}
\label{eq:ball}
|A(B) - V(B)| ~=~ k^{O(1)} \cdot V_{k-3}(t)~.
\end{equation}
We prove (\ref{eq:ball}) in the sequel.

Another subtle point is that we have rotated the vector space to move
from the standard Kronecker basis to the orthonormal basis $\Upsilon$.
(In fact, $\Upsilon$ is an orthonormal basis for the hyperplane $\cH'$,
but it can be completed to an orthonormal basis for $\MathR^k$ by
inserting the vector ${\delta \over {||\delta||}}$ into it.)
Consequently, the integer lattice was rotated as well, and so in our
context $A(B)$ is actually the number of points of the {\em rotated }
integer lattice that are contained in $B$.
These two quantities may be slightly different. However, we argue
below that the estimate (\ref{eq:ball}) applies for the rotated
integer lattice as well, for any rotation.

Recall that $m = |I|$.
Let 
\begin{eqnarray}
\label{eq:Qext}
\tQ_{\mathit{ext}} & = & \{\alpha =
(\alpha_1,\alpha_2,\ldots,\alpha_{k-1}) \in \MathR^{k-1} :
|| \alpha||^2 \le T', ~\alpha_i \ge 0 \mbox{~for~all~} i \ge m\} \\
\label{eq:Qint}
\tQ_{\mathit{int}} & = & \{\alpha =
(\alpha_1,\alpha_2,\ldots,\alpha_{k-1}) \in \MathR^{k-1} : \\
\nonumber
&& 
|| \alpha||^2 \le T' - (g + 1), ~\alpha_i \ge 0 \mbox{~for~all~} i \ge m\}
\end{eqnarray}
Observe that $\tQ \subseteq \tQ_{\mathit{ext}} \setmns
\tQ_{\mathit{int}}$.
Also, let $\tZ$ denote 
\begin{equation}
\label{eq:tZ}
\tZ ~=~ \{ \alpha = (\alpha_1,\alpha_2,\ldots,\alpha_{k-3}) \in
\MathR^{k-3} : ||\alpha||^2 \le T',  \mbox{~for~all~} i \ge m\}~.
\end{equation}
The set $ \tQ_{\mathit{ext}}$ (respectively, $\tQ_{\mathit{int}}$) is
the intersection of the $(k-1)$-dimensional ball of squared radius
$T'$ (resp., $T' - (g+1)$) centered at the origin with the half-spaces
$\alpha_i \ge 0$ for all $i \ge m$. The set $\tZ$ is the intersection
of the $(k-3)$-dimensional ball  of squared radius
$T'$ centered at the origin with the half-spaces
$\alpha_i \ge 0$ for all $i \ge m$.
The analogue of (\ref{eq:ball}) that is required for our argument is
\begin{equation}
\label{eq:cut_ball}
|A(\tQ_{\mathit{ext}} ) - \Vl(\tQ_{\mathit{ext}} )| ~=~
 k^{O(1)} \cdot \Vl(\tZ)~.
\end{equation}

Given (\ref{eq:cut_ball}) we show (\ref{eq:Q}) by the following
argument.

\begin{lemma}
\label{lm:AQ}
$A(\tQ) = 2^{O(\sqrt{k})} \cdot 2^{\eps k} \cdot \left({{\pi e} \over
  6}\right)^{k/2} \cdot y^{k-3}$.
\end{lemma}
\proof
By (\ref{eq:cut_ball}),
\begin{eqnarray*}
A(\tQ)  & \le & A(\tQ_{\mathit{ext}}) -  A(\tQ_{\mathit{int}}) ~\le~
\Vl(\tQ_{\mathit{ext}}) + k^{O(1)} \cdot \Vl(\tZ) -
\Vl(\tQ_{\mathit{int}}) 
+ k^{O(1)} \cdot \Vl(\tZ) \\
& = & ( \Vl(\tQ_{\mathit{ext}}) - \Vl(\tQ_{\mathit{int}}) )
+ k^{O(1)} \cdot \Vl(\tZ).
\end{eqnarray*}
Observe that 
$\Vl(\tZ) = {{2^{\eps k}} \over {2^{k-3}}} \cdot \beta_{k-3} \cdot
  (T')^{{k-3} \over 2}$.
Also, since $T'$ is much greater than $g$, 
\begin{eqnarray*}
&& \Vl(\tQ_{\mathit{ext}}) - \Vl(\tQ_{\mathit{int}}) ~=~
{{2^{\eps k}} \over {2^{k-1}}} \cdot \beta_{k-1} ((T')^{{k-1} \over 2} - 
(T' - (g+1))^{{k-1} \over 2}) \\
&& ~\le ~ O(1) \cdot {{2^{\eps k}} \over {2^{k-1}}} \cdot \beta_{k-1} \cdot
k \cdot (g+1) \cdot (T')^{{k-3} \over 2}~.
\end{eqnarray*}
Hence
$$ 
A(\tQ) ~\le~ O(1) \cdot {{2^{\eps k}} \over {2^{k-3}}} \cdot  (k^{O(1)} \cdot
\beta_{k-3} + k \cdot (g+1) \cdot \beta_{k-1}) \cdot (T')^{{k-3} \over 2}~.
$$
Since $\beta_{k-3} = \Theta(k \cdot \beta_{k-1})$ and $g \le k$, it
follows that 
$$A(\tQ) ~\le~ k^{O(1)} \cdot {{2^{\eps k}} \over {2^{k-1}}} \cdot 
\beta_{k-1} 
\cdot (T')^{{k-3} \over 2}~.
$$
By (\ref{eq:ubR}), $T' \le {k \over 3} \cdot y^2 (1 + O({1 \over
  {\sqrt{k}}}))$.
Also, $\beta_{k-1} = {{\pi^{{k-1} \over 2}} \over \Gamma({{k+1} \over
    2})}$.
Hence
$$A(\tQ) = 2^{O(\sqrt{k})} \cdot 2^{\eps k} \cdot \left({{\pi e} \over
  6}\right)^{k/2} \cdot y^{k-3}~.
\inQED
$$

Hence it remains to prove (\ref{eq:cut_ball}).
Our proof is closely related to the argument in \cite{F82}, pp. 94-97,
 and is provided for the sake of completeness.
In addition, our argument is more general than the one in \cite{F82},
 as the latter argument 
applies only for balls, while our argument applies for intersections
of balls with half-spaces.


Fix $m$ to be a positive integer number. (In our application $m =
|I|$.)
For positive integer numbers $k$ and $t$, let $Q_k(t)$ denote the
intersection of the $k$-dimensional ball $B_k(t)$ centered at the
origin with squared radius $t$ with the half-spaces $\cH^{(i)} = \{\alpha =
(\alpha_1,\alpha_2,\ldots,\alpha_k) \mid \alpha_i \ge 0\}$, 
for all $i \ge m$.
Let $\bV_k(t)$ denote the volume $\Vl(Q_k(t))$, and $\bA_k(t)$ denote
the number of points of the rotated integer lattice in $Q_k(t)$.
Note that 
$\bV_k(t) = {{\beta_k} \over {2^{\max \{k - m + 1,0\}}}} \cdot t^{k/2}$.
The next lemma provides an upper bound for the discrepancy between
$\bV_k(t)$ and $\bA_k(t)$ in terms of $\bV_{k-2}(t)$.

\begin{lemma}
\label{lm:discrep}
For a sufficiently large real $t > 0$ and an integer $k \ge 5$, 
%
$$|\bA_k(t) - \bV_k(t)| = O(k^{3/2} \cdot \bV_{k-2}(t))~.$$
\end{lemma}
\inline Remark: This lemma applies even if $k = k(t)$ is a function of $t$.

Before proving Lemma \ref{lm:discrep}, we first provide a number of
auxiliary lemmas that will be useful for its proof.
We start with Euler Sum-formula (\cite{F82}, Satz 29.1, p.185).

\begin{lemma} 
\label{lm:EulerSum}
For a real-valued function $f(u)$ differentiable in a segment $[a,b]$,
$$
\sum_{a < \ell \le b} f(\ell) ~=~ \int_a^b f(u) du + \psi(a) \cdot
f(a)  - \psi(b) \cdot f(b) + \int_a^b \psi(u) \cdot f'(u) du~,$$
where $\psi(u) = u - \lfloor u \rfloor - \half$.
\end{lemma}

In addition, we will use the following property of the function
$\psi(\cdot)$.

\begin{lemma}
\label{lm:psi}
For any two real numbers $\kappa$ and $\lambda$, $\kappa \le \lambda$,
$- \half ~\le~ \int_\kappa^\lambda \psi(u) du ~\le ~ 1$.
\end{lemma}
\proof
Observe that
$$\int_0^1 \psi(u) du ~=~ \int_0^1 \left(u - \lfloor u \rfloor  - \half\right)
du 
~=~
\int_0^1 \left(u - \half \right)du ~=~ \int_{-\half}^{\half} t dt ~=~ 0~.$$
Moreover, for any integer $j$,
$\int_j^{j+1} \psi(u) du ~=~ \int_j^{j+1} \left(u - \lfloor u \rfloor -
\half \right) du ~=~  \int_{-\half}^{\half} t dt ~=~ 0$.
Hence $\int_0^j \psi(u) du = \int_{-j}^0 \psi(u) du = 0$, for any
positive integer $j$.

It follows that
$$\int_\kappa^\lambda \psi(u) du ~=~ \int_\kappa^{\lceil \kappa
  \rceil} \psi(u) du + \int_{\lceil \kappa \rceil}^{\lfloor \lambda
  \rfloor} \psi(u) du + \int_{\lfloor \lambda
  \rfloor}^\lambda \psi(u) du ~=~ \int_\kappa^{\lceil \kappa
  \rceil} \psi(u) du + \int_{\lfloor \lambda
  \rfloor}^\lambda \psi(u) du~.
$$
Let $\kappa' = \kappa - \lfloor \kappa \rfloor$, $\lambda' = \lambda -
\lfloor \lambda \rfloor$.
The integral $\int_\kappa^{\lceil \kappa
  \rceil} \psi(u) du$ is equal to 0 if $\kappa = {\lceil \kappa
  \rceil}$, and is equal to $\int_{\kappa'}^1 \psi(u) du$, otherwise.
The latter integral satisfies
$$
\int_{\kappa'}^1 \psi(u) du ~=~ \int_{\kappa' - \half}^{\half} t dt
~=~ \half (1 - \kappa') \kappa' ~\le~ \half~.
$$
Also, $\int_{\kappa'}^1 \psi(u) du = \half (1 - \kappa') \kappa'  \ge 0$.

Analogously, $\int_{\lfloor \lambda
  \rfloor}^\lambda \psi(u) du$ is equal to 0 if $\lambda = \lfloor
  \lambda \rfloor$, and it is equal to $\int_0^{\lambda'} \psi(u) du$,
otherwise.
The latter integral satisfies
$$
\int_0^{\lambda'} \psi(u) du ~=~ \int_{-\half}^{\lambda' - \half} t dt
~=~
\half (\lambda' - 1)\lambda~.$$
Hence
$-\half \le \int_0^{\lambda'} \psi(u) du \le \half$,
and thus,
$-\half \le \int_\kappa^\lambda \psi(u) du \le 1$.
\QED

Next, we use Lemma \ref{lm:psi}   to derive another useful property of
the function $\psi(\cdot)$.

\begin{lemma}
\label{lm:psi_integral}
For a positive real number $t$ and a positive integer $p \ge 2$,
$$|\int_0^{\sqrt{t}} u \cdot \psi(u) (t -u^2)^{{p \over 2} -1} du| ~\le~
t^{{p-1} \over 2}~.$$
\end{lemma}
\proof
Since $f(u) = u$ is a monotone increasing function, there exists $\xi
\in [0,\sqrt{t}]$ such that 
$$
\int_0^{\sqrt{t}} u \cdot \psi(u) (t -u^2)^{{p \over 2} -1} du ~=~
\sqrt{t} \int_\xi^{\sqrt{t}} \psi(u)  (t -u^2)^{{p \over 2} -1} du~.
$$
Since $g(u) = (t -u^2)^{{p \over 2} -1} $ is a monotone decreasing
function in $[\xi,\sqrt{t}]$, 
there exists $\eta \in [\xi,\sqrt{t}]$
such that the right-hand-side  is equal
to $\sqrt{t} \cdot (t - \xi^2)^{{p \over 2} - 1} \int_\xi^\eta \psi(u) du$.
By Lemma  \ref{lm:psi}, 
$$|\sqrt{t} (t - \xi^2)^{{p \over 2} - 1} \int_\xi^\eta \psi(u) du|
~\le~
\sqrt{t} \cdot t^{{p \over 2}  -1} ~=~ t^{{p-1} \over 2}~.
\inQED
$$
Note also that
$$
|\int_{-\half}^0  u \cdot \psi(u) (t -u^2)^{{p \over 2} -1} du|
~\le~
|\int_{-\half}^0  (t -u^2)^{{p \over 2} -1} du| ~\le~ \half \cdot t^{{p
    \over 2} - 1}~.
$$
Hence for any $t$ and $p$ as above,
\begin{equation}
\label{eq:integral}
|\int_{-\half}^{\sqrt{t}}  u \cdot \psi(u) (t -u^2)^{{p \over 2} -1} du|
~\le~ t^{{p-1} \over 2} + \half \cdot t^{{p \over 2} - 1} ~\le~
t^{{p-1} \over 2} \left(1 + {1 \over {2 \sqrt{t}}}\right)~.
\end{equation}
We are now ready to prove Lemma \ref{lm:discrep}.

\inline Proof of Lemma \ref{lm:discrep}:

We prove by induction on $k$ that there exists a universal constant
$c > 0$ such that 
\begin{equation}
\label{eq:ind_hyp}
|\bA_k(t) - \bV_k(t) | \le \left(c \cdot \sum_{j=1}^{k-1}
\sqrt{j} \right)\left( 1 + {1 \over
  {2\sqrt{t}}}\right) \cdot \bV_{k-2}(t)~.
\end{equation}
The constant $c$ will be determined later.

The induction base is $k=5$.
It is well-known (see, e.g., \cite{Adh94}) 
that $|\bA_5(t) - \bV_5(t)| = O(\bV_3(t)) =
O(t^{3/2})$.

Next, we prove the induction step.

In all summations below, $\ell$ is an integer index.
The analysis splits into two cases.
The first case is $k + 1 < m$, and the second is $k+1 \ge m$.
In the first case
$$
\bA_{k+1}(t) ~=~ \sum_{|\ell| \le \sqrt{t}} \bA_k(t - \ell^2) ~=~
\sum_{|\ell| \le \sqrt{t}} \bV_k(t - \ell^2) + 
\sum_{|\ell| \le \sqrt{t}} (\bA_k(t - \ell^2) - \bV_k(t - \ell^2))~,
$$
and so
$$
|\bA_{k+1}(t) - \sum_{|\ell| \le \sqrt{t}} \bV_k(t - \ell^2)| ~=~
|\sum_{|\ell| \le \sqrt{t}} (\bA_k(t - \ell^2) - \bV_k(t - \ell^2))|
~\le~ \sum_{|\ell| \le \sqrt{t}} |\bA_k(t - \ell^2) - \bV_k(t -
\ell^2)|~.
$$
In the second case the same inequalities apply, but the index $\ell$
runs in the range $0 \le \ell \le \sqrt{t}$ in all summations.
It turns out to be more convenient to have the index $\ell$ vary in
the range
 $-\half \le \ell \le \sqrt{t}$ rather than $0 \le \ell \le \sqrt{t}$
in these summations.

By the induction hypothesis (that is, by (\ref{eq:ind_hyp})),
$$|\bA_k(t - \ell^2) - \bV_k(t - \ell^2) | ~\le~
\left(c \cdot \sum_{j=1}^{k-1}
\sqrt{j} \right)\left( 1 + {1 \over
  {2\sqrt{t}}}\right) \cdot \bV_{k-2}(t - \ell^2)~.$$
Hence
\begin{equation}
\label{eq:A_est}
|\bA_{k+1}(t) - \sum_{|\ell| \le \sqrt{t}} \bV_k(t - \ell^2)| 
~\le~
\left(c \cdot \sum_{j=1}^{k-1}
\sqrt{j} \right)\left( 1 + {1 \over
  {2\sqrt{t}}}\right)  \cdot \sum_{|\ell| \le
  \sqrt{t}} \bV_{k-2}(t - \ell^2)~.
\end{equation}
Next, we estimate  $\sum_{|\ell| \le \sqrt{t}} \bV_k(t - \ell^2)$
via Euler Sum-formula (Lemma \ref{lm:EulerSum}).
In the first case,
we substitute $a = -\sqrt{t}$, $b = \sqrt{t}$, and $f(u) = \bV_k(t -
u^2)$.
Then $f(a) = f(b) = \bV_k(0) = 0$,
and 
$${{df} \over {du}}(u) ~=~ {d \over {du}} \bV_k(t - u^2) ~=~
\beta_k {d \over {du}} (t - u^2)^{k \over 2} ~=~
- \beta_k \cdot k \cdot  (t - u^2)^{{k \over 2} - 1} u~.$$
By Lemma \ref{lm:EulerSum} it follows that
\begin{equation}
\label{eq:sumV_fst}
\sum_{|\ell| \le \sqrt{t}} \bV_k(t - \ell^2) ~=~
\int_{\sqrt{t}}^{\sqrt{t}} \bV_k(t - u^2) du - k \cdot \beta_k
\int_{\sqrt{t}}^{\sqrt{t}} \psi(u) (t - u^2)^{{k \over 2} - 1} u 
du~.
\end{equation}
In the second case ($k \ge m-1$),
$a = - \half$, $b = \sqrt{t}$, and again $f(a) = f(b) = 0$.
Also,
$${{df} \over {du}}(u) ~=~ - \beta_k \cdot {1 \over {2^{k -m + 1}}}
\cdot k
\cdot (t - u^2)^{{k \over 2} - 1} u~,$$
and thus,
\begin{equation}
\label{eq:sumV_scnd}
\sum_{-\half < \ell \le \sqrt{t}} \bV_k(t - \ell^2) ~=~
\int_{-\half}^{\sqrt{t}} \bV_k(t - u^2) du - k \cdot {{\beta_k} \over  
{2^{k -m + 1}}}
\int_{-\half}^{\sqrt{t}} \psi(u) (t - u^2)^{{k \over 2} - 1} u  du~.
\end{equation}
In the first case, since $h(u) = u \psi(u)$ is an even function on $\MathR
\setmns \MathZ$, the right-hand-side in (\ref{eq:sumV_fst}) is equal
to
$$
\int_{\sqrt{t}}^{\sqrt{t}} \bV_k(t - u^2) du - 2k \cdot \beta_k
\int_{0}^{\sqrt{t}} \psi(u) (t - u^2)^{{k \over 2} - 1} u  du~.
$$
Let $J$ denote $|\int_0^{\sqrt{t}} u \cdot \psi(u) (t - u^2)^{{k \over
    2} - 1} du |$.
By Lemma \ref{lm:psi_integral},
$
J
\le
t^{{k-1} \over 2}
$.
Hence
$$
\sum_{|\ell| \le \sqrt{t}} \bV_k(t - \ell^2) ~=~ 
\int_{\sqrt{t}}^{\sqrt{t}} \bV_k(t - u^2) du - 2 k \beta_k \cdot J ~=~
\bV_{k+1}(t) - 2 k \beta_k \cdot J~.
$$
It follows that
\begin{equation}
\label{eq:VV_fst}
|\bV_{k+1}(t) - \sum_{|\ell| \le \sqrt{t}} \bV_k(t - \ell^2)| ~=~ 2 k
\beta_k \cdot J ~\le~ 2k \beta_k \cdot t^{{k-1} \over 2}~.
\end{equation}
In the second case by  (\ref{eq:sumV_scnd}) and since
$\int_{0}^{\sqrt{t}} \bV_k(t - u^2) du = \bV_{k+1}(t)$, it follows that
$$
\sum_{-\half < \ell \le \sqrt{t}} \bV_k(t - \ell^2) ~ = ~
\int_{-\half}^0 \bV_k(t - u^2) du + \bV_{k+1}(t) - k \cdot {{\beta_k} \over
  {2^{k - m + 1}}} \int_{-\half}^{\sqrt{t}} \psi(u) u (t - u^2)^{{k
    \over 2} - 1} du~.
$$
Let $J'$ denote $| \int_{-\half}^{\sqrt{t}}  u \cdot \psi(u) (t -u^2)^{{k \over 2} -1} du|$.
By (\ref{eq:integral}), 
\begin{equation}
\label{eq:Jtag}
J' 
~\le~ t^{{k-1} \over 2} + \half \cdot t^{{k \over 2} - 1} ~\le~
t^{{k-1} \over 2} \left(1 + {1 \over {2 \sqrt{t}}}\right)~.
\end{equation}
Since $\bV_k(t - u^2) \ge 0$ for all $u$, $-\half \le u \le
0$,
the integral $\int_{-\half}^0 \bV_k(t - u^2) du$ is non-negative as well.
Thus,
$$
|\bV_{k+1}(t) - \sum_{-\half \le \ell \le \sqrt{t}} \bV_k(t - \ell^2)| ~\le~  k
 \cdot {{\beta_k} \over {2^{k - m +1}}} \cdot J' 
~\le~ k \cdot {{\beta_k} \over {2^{k - m +1}}} 
\cdot t^{{k-1} \over 2} \left(1 +
 {1 \over {2 \sqrt{t}}} \right)~.
$$
In the first case, 
by the triangle inequality, by (\ref{eq:A_est}), and by
 (\ref{eq:VV_fst}),
\begin{eqnarray}
\label{eq:triang}
&& |\bA_{k+1}(t) - \bV_{k+1}(t)| ~\le ~
|\bA_{k+1}(t) - \sum_{|\ell| \le \sqrt{t}} \bV_k(t - \ell^2)| 
\\
\nonumber
&+ & | \sum_{|\ell| \le \sqrt{t}} \bV_k(t - \ell^2) - \bV_{k+1}(t) | \\
\nonumber
& \le & \left(c \cdot \sum_{j=1}^{k-1}
\sqrt{j} \right)
 \sum_{|\ell| \le
  \sqrt{t}} \bV_{k-2}(t - \ell^2) 
+ | \sum_{|\ell| \le \sqrt{t}} \bV_k(t - \ell^2) - \bV_{k+1}(t) | \\
\label{eq:bnd}
& \le & \left(c \cdot \sum_{j=1}^{k-1}
\sqrt{j} \right)
 \sum_{|\ell| \le
  \sqrt{t}} \bV_{k-2}(t - \ell^2) + 2 k  \beta_k \cdot t^{{k-1}
  \over 2}~. 
\end{eqnarray}

Analogously, in the second case,
\begin{eqnarray}
\label{eq:triang_scnd}
|\bA_{k+1}(t) - \bV_{k+1}(t)| &\le & 
\left(c \cdot \sum_{j=1}^{k-1}
\sqrt{j} \right) \left( 1 + {1 \over
 {2\sqrt{t}}}\right) \sum_{-\half \le \ell \le \sqrt{t}}
\bV_{k-2}(t - \ell^2)  \\
\nonumber
& + & 2 k \cdot {{\beta_k} \over {2^{k + 2 - m}}}
\cdot t^{{k-1} \over 2} \cdot \left(1 +  {1 \over
 {2\sqrt{t}}}\right)~.
\end{eqnarray}
However, in the first case 
$\sum_{|\ell| \le \sqrt{t}} \bV_{k-2}(t - \ell^2) ~\le~
\int_{-\sqrt{t}}^{\sqrt{t}} \bV_{k-2}(u) du ~=~ \bV_{k-1}(t)$.
In the second case,
$$
\sum_{-\half \le \ell \le \sqrt{t}} \bV_{k-2}(t - \ell^2) ~=~
\sum_{0  \le \ell \le \sqrt{t}} \bV_{k-2}(t - \ell^2) ~\le~
\int_0^{\sqrt{t}} \bV_{k-2}(u) du ~=~ \bV_{k-1}(t)~.
$$
Hence in both cases the first terms in (\ref{eq:bnd}) and in the
right-hand-side of (\ref{eq:triang_scnd}) are at most
$$
\left(c \cdot \sum_{j=1}^{k-1}
\sqrt{j} \right) \left( 1 + {1 \over
 {2\sqrt{t}}}\right) \bV_{k-1}(t)~.
$$
Consequently,  in both cases,
\begin{eqnarray*}
|\bA_{k+1}(t) - \bV_{k+1}(t) | &\le &
\left(c \cdot \sum_{j=1}^{k-1}
\sqrt{j} \right) \left( 1 + {1 \over
 {2\sqrt{t}}}\right) \cdot \bV_{k-1}(t) \\
& + & 2k \cdot  {{\beta_k} \over {2^{\max\{k + 2 - m,0\}}}}
\cdot t^{{k-1} \over 2} \cdot \left(1 +  {1 \over
 {2\sqrt{t}}}\right)~.
\end{eqnarray*}
By (\ref{eq:def_beta}), 
 $\beta_k = \Theta\left({{\beta_{k-1}} \over {\sqrt{k}}}\right)$.
Set
 $c$ to be a universal constant such that $c \ge {{{\sqrt{k} \cdot
      \beta_k}} \over {2 \beta_{k-1}}}$, for all integer $k \ge 2$.
Then
\begin{eqnarray*}
&& |\bA_{k+1}(t) - \bV_{k+1}(t) | ~\le  ~ 
\left(c \cdot \sum_{j=1}^{k-1}
\sqrt{j} \right) \left( 1 + {1 \over
 {2\sqrt{t}}}\right) \cdot \bV_{k-1}(t) 
\\
& + &
 c
\cdot \sqrt{k} \cdot \left( 1 + {1 \over
 {2\sqrt{t}}}\right) \cdot {{\beta_{k-1}} \over {2^{\max\{(k - 1) - m
      + 1, 0\}}}} \cdot t^{{k-1} \over 2} 
~ =~  \left(c \cdot \sum_{j=1}^{k-1}
\sqrt{j} \right)  \left( 1 + {1 \over
 {2\sqrt{t}}}\right) \bV_{k-1}(t)~.
\end{eqnarray*}
Finally, $\sum_{j = 1}^k \sqrt{j} \le k^{3/2}$, completing the proof.
\QED


\section{Conclusion}
\label{sec:concl}

In this paper we improved the lower bound of Behrend by a factor of
$\Theta(\sqrt{\log n})$. As was already mentioned, both Behrend's and
our proof arguments rely on the Pigeonhole Principle.
It is reasonable to believe that by choosing $T = R^2 = \mu_Z$ (see 
(\ref{eq:Rdef})) one can get an annulus with at least as many integer
points as in the annulus $\cS$ chosen via the Pigeonhole Principle. To
prove that this is the case one should probably use normal
approximation of the discrete random variable $Z$ (see Sections 
\ref{sec:behrend}
and \ref{sec:our}), and employ probablistic estimates to argue that
the probability that $Z$ is between $(\mu_Z - {{\eps k} \over 2})$ 
and $(\mu_Z + {{\eps k} \over 2})$ is at least as large as the
probability that it is between $(\mu_Z - 2\sigma_Z)$ and $(\mu_Z + 2
\sigma_Z)$, divided by ${{\eps k} \over {4 \sigma_Z}}$. Although this
appears to be quite clear intuitively, so far we were not able to find
sufficiently precise probabilistic estimates to prove this statement
formally.
Once this intuition is formalized, our construction will become
independent of the Pigeonhole Principle. This, in turn, would be a
significant improvement of the lower bound of Moser \cite{M53}.

\section*{Acknowledgements}

The author is indebted to Don Coppersmith, who was offered a
coauthorship on this paper. In particular, Lemma \ref{lm:wipe} is due
to Don. In addition, fingerprints of Don can be found in numerous
other places in this paper. 

The author is grateful to Benny Sudakov for introducing him to the
problem.
The author thanks also Noga Alon, Eitan Bachmat,
B\'ela Bollob\'as,  Danny Berend, Alexander Razborov, Oded Regev,
Alex Samorodnitsky, and Shakhar Smorodinsky, for encouragement and for
helpful discussions.

\bibliographystyle{abbrv}{
\bibliography{progfree}  
}


\end{document}